\documentclass[12pt]{amsart}

\newtheorem{satz}{Theorem}[section]

\newtheorem{lemma}[satz]{Lemma}
\newtheorem{theo}[satz]{Theorem}

\newtheorem{prop}[satz]{Proposition}

\newtheorem{cor}[satz]{Corollary}

\newtheorem{df}[satz]{Definition}

\begin{document}


\date{}
\title[Separable Sets in Universal Algebra]{\bf Essential Variables and Separable Sets in Universal Algebra}

\author{Sl. Shtrakov}
\address{Sl. Shtrakov\\
South-West-University Blagoevgrad\\
Faculty of Mathematics and Sciences\\
Blagoevgrad\\
Bulgaria}
\email{\ttfamily shtrakov@aix.swu.bg}
\urladdr{ http://home.swu.bg/shtrakov/}
\author{K. Denecke}
\address{K. Denecke\\
Universit\"at Potsdam\\
Fachbereich Mathematik\\
Postfach 601553\\
D-14415 Potsdam}
\email{\ttfamily kdenecke@rz.uni-potsdam.de}
\urladdr{ http://users.math.uni-potsdam.de/\~denecke/}

\maketitle

\begin{abstract}
The study of essential and strongly essential variables in functions defined
on finite sets is a part of  $k$-valued logic. We extend the main definitions
from functions to terms. This allows us to apply concepts and results of
Universal Algebra. On the basis of the concept of a separable set of
variables in a term we introduce a new notion of complexity of terms,
algebras and varieties and give examples.

\end{abstract}

\section{Introduction}\label{sec1}
\noindent

Several authors considered essential  variables of
functions under different aspects ({\it see e.g. \cite{s1,s2,s3,s5}}).
The unary
function $f: A \rightarrow A$ {\it depends essentially} on its input $x$ if
it takes on at least two values, i.e. if $f$ is not constant. The $n$-ary
function $f: A^n \rightarrow A$ depends essentially on its $i$-the input
$x_i$ if there are elements $a_1, \ldots,a_{i-1}, a_{i+1}, \ldots, a_n \in
A$ such that the unary function defined by
\[x_i \longmapsto f(a_1, \ldots, a_{i-1}, x_i, a_{i+1}, \ldots a_n) \]
is not constant on $A$.

In this case it is very common to say, the variable $x_i$ is essential for
$f(x_1, \ldots, x_i, \ldots, x_n)$ and to consider the set $Ess(f)$ of all
variables which are essential for $f(x_1, \ldots, x_n)$. As
usual, instead of the function $f$ one writes the term $f(x_1, \ldots,  x_n)$.
So, it is very natural to define essential variables for
terms. One can do this by using the former definition and going from terms
to term operations which are induced by terms on a given algebra
$\mathcal U$. This is a new point of view and allows us to apply the
methods of Universal Algebra to the study of essential variables.

We will give the necessary definitions and prove some consequences.
Moreover, on the basis of essential variables we define a new concept of
complexity for terms, for polynomials, and for algebras. This concept
depends not only on the
syntax of the term but also on its meaning in a given model.

\section{Basic concepts}\label{sec2}

\noindent
We will use the denotation  $O_A^n$ for the set of all $n$-ary functions
defined on the set $A$, i.e.,  $O_A^n := \{f~|~f: A^n \rightarrow A \}$
and $O_A := \bigcup\limits_{n=1}^{\infty} O_A^n$.
An {\it algebra of type $\tau$} is
a pair ${\mathcal U} = (A; (f_i^{\mathcal U})_{i \in I})$ where $(f_i^{\mathcal U})_{i \in I}$
is an indexed subset of $O_A$ and $f_i^{\mathcal U}$ is $n_i$-ary, $n_i \geq 1$. By
$Alg(\tau)$ we denote the class of all algebras of type $\tau$.

{\it Terms of type $\tau$} are defined in the following inductive way. Let
$(f_i)_{i \in I}$ be an indexed set of operation symbols where $f_i$ is
$n_i$-ary and let $X_n = \{ x_1, \ldots, x_n \}$ be an $n$-element set of
variables. Then we define {\it $n$-ary terms of type $\tau$} as follows
\begin{enumerate}
\item[(i)] $x_i \in X_n$ is an $n$-ary term for all $i \in \{ 1, \ldots, n \}$,
\item[(ii)] if $t_1, \ldots, t_{n_i}$ are $n$-ary terms and if $f_i$ is $n_i$-ary
      then $f_i(t_1, \ldots, t_{n_i})$ is an $n$-ary term.\\
\end{enumerate}

By $W_\tau(X_n)$ we denote the set of all $n$-ary terms of type $\tau$. If
$X = \{ x_1, \ldots, x_n, \ldots \}$ is a countably infinite alphabet then
$W_\tau(X) := \bigcup\limits_{n=1}^{\infty} W_\tau(X_n)$ is the set of all
terms of type $\tau$. For $t \in W_\tau(X)$ by $var(t)$ we denote the set of
all variables which occur in $t$.

Polynomials of type $\tau$ are defined in a similar way. Besides the set $X$
we need a set ${\mathcal A}$ of constant symbols with ${\mathcal A} \cap X
= \emptyset$. {\it Polynomials of type $\tau$ over ${\mathcal A}$} are
defined in three steps:
\begin{enumerate}
\item[(i)] variables from $X$ are polynomials of type $\tau$ over ${\mathcal A}$,
\item[(ii)] elements ${\overline a}$ from ${\mathcal A}$ are polynomials of type
      $\tau$ over ${\mathcal A}$,
\item[(iii)] if $p_1, \ldots, p_{n_i}$ are polynomials of type $\tau$ over
      ${\mathcal A}$ and if $f_i$ is an $n_i$-ary operation symbol then
      $f_i(p_1, \ldots, p_{n_i})$ is a polynomial of type $\tau$ over
      ${\mathcal A}$.
\end{enumerate}
Let $P_\tau(X, {\mathcal A})$ be the set of all polynomials of type $\tau$
over ${\mathcal A}$. Clearly, $W_\tau(X) \subset P_\tau(X, {\mathcal A})$.
Let ${\mathcal B}$ be an algebra of type $\tau$ containing a subalgebra
${\mathcal U}$ whose universe has the same cardinality as ${\mathcal A}$.
Then every $n$-ary polynomial $p$ of type $\tau$ over ${\mathcal A}$
induces on the algebra ${\mathcal B}$ an $n$-ary polynomial operation
$p^{\mathcal B}$ defined by the following steps:
\begin{enumerate}
\item[(i)] if $x_i \in X_n$,  then $x_i^{\mathcal B} := e_i^{n,\mathcal B},~1 \leq i \leq n$, where
      $e_i^{n,\mathcal B}: (b_1, \ldots, b_n)$ $\longmapsto b_i$ is the $n$-ary
      projection onto the $i$-th coordinate,
\item[(ii)] if ${\overline a} \in {\mathcal A}$ then ${\overline a}^{\mathcal B} := c_a^n$ is the $n$-ary constant operation on $\mathcal B$
      with value $a \in A$
      and every element from $A \subseteq B$ is uniquely induced by an
      element from ${\mathcal A}$,
\item[(iii)] if $p = f_i(p_1, \ldots, p_{n_i})$ and if we assume that the
      polynomial operations $p_1^{\mathcal B}, \ldots, p_{n_i}^{\mathcal B}$ are already defined,
      then $p^{\mathcal B} = f_i^{\mathcal B}(p_1^{\mathcal B},
      \ldots, p_{n_i}^{\mathcal B})$ where the right hand side is the
      usual composition of operations defined by
      \[f_i^{\mathcal B}(p_1^{\mathcal B}, \ldots, p_n^{\mathcal
        B})(b_1, \ldots, b_n) := f_i^{\mathcal B}(p_1^{\mathcal B}(b_1,
        \ldots, b_n), \ldots, p_{n_i}^{\mathcal B}(b_1, \ldots, b_n)). \]
\end{enumerate}
The set of all operations induced by arbitrary polynomials of type $\tau$
over ${\mathcal A}$ on the algebra ${\mathcal B}$ is denoted by
$P_A({\mathcal B})$. By $T({\mathcal B})$ we denote the set of all
operations induced by arbitrary terms from $W_\tau(X) \subset P_\tau(X,
{\mathcal A})$. The elements from $T({\mathcal B})$ are called {\it term
operations} of ${\mathcal B}$ and $T({\mathcal B})$ is called {\it clone
of~ term operations} of ${\mathcal B}$. The elements of $P_A({\mathcal B})$ are called polynomial operations of ${\mathcal B}$. Polynomial
operations can also be defined by mappings $h: X \cup {\mathcal A}
\rightarrow B$ which are extensions of a distinguished mapping
$h': \mathcal A \to A$ with
$h'({\overline a}) = a \in A \subseteq B$ for ${\overline a} \in {\mathcal A}$.

Any such mapping $h$ is called evaluation of $X \cup {\mathcal A}$ with $B$. It
is well-known that an evaluation mapping can be uniquely extended to a mapping
\[{\overline h}: P_\tau(X, {\mathcal A}) \rightarrow  B~.\]
Here for operation symbols $f_i$ occurring in a polynomial $p$ one has to
substitute the corresponding operations $f_i^{\mathcal B}$.
For terms the extension ${\overline h}$ maps $W_\tau(X)$ to $\mathcal B$.

If $s,t$ are terms of type $\tau$ then $s \approx t$ is called identity
satisfied in the algebra ${\mathcal B}$ of type $\tau$ if the induced term
operations are equal, i.e. if $s^{\mathcal B} = t^{\mathcal B}$.
In this case we write $\mathcal{B} \models s \approx t$. Polynomial identities in $\mathcal B$ are pairs of polynomials
$p\approx q$ such that the polynomial operations induced by $p$ and $q$ on $\mathcal B$ are equal.
A variety $V$ of type $\tau$ is a class of algebras of type $\tau$ such that
there exists a set $\Sigma$ of equations of type $\tau$ with the property
that $V$ consists exactly of all algebras of type $\tau$ such that every
equation from $\Sigma$ is satisfied as identity.

As usual, by $\Bbb{H}, \Bbb{S}, \Bbb{P}$ we denote the operators of forming
homomorphic images, subalgebras, and direct products of a given algebra or a
given class of algebras. A class $V$ of algebras of the same type $\tau$ is a
variety iff $V$ is closed under the operators $\Bbb{H}, \Bbb{S}$, and
$\Bbb{P}$.


\section{Essential variables and separable sets of variables in terms with
respect to an algebra}

\begin{df}\label{3.1}\rm Let $t \in W_\tau(X_n)$ be an $n$-ary term of type
$\tau$ and let $\mathcal U$ be an algebra of type $\tau$. Then the
variable $x_i,~1 \leq i \leq n$, is called {\it essential in $t$ with
respect to the algebra ${\mathcal U}$} if the term operation
$t^{\mathcal U}: A^n \rightarrow A$ induced by $t$ on the algebra
${\mathcal U}$ depends essentially on its $i$-th input $x_i$.
By $Ess(t, {\mathcal U})$ we denote the set of all variables which are
essential in $t$ with respect to the algebra ${\mathcal U}$. For a
polynomial $p$ we define in the same way when $x_i$ is essential in $p$ with
respect to ${\mathcal U}$.
\end{df}
\noindent
{\it Remarks}

\noindent
\begin{enumerate}
\item The definition means that $x_i$ is essential in the term $t$ with
      respect to
      ${\mathcal U}$ if there are two different evaluation mappings $h,h':
      X_n \rightarrow A$ with $h/X_n \setminus \{x_i \} =
      h'/X_n \setminus \{ x_i \}$ such that ${\overline h}(t) \not=
      {\overline h'}(t)$ where ${\overline h}: W_\tau(X_n) \rightarrow
       A$ is the extension of $h$.
\item If ${\mathcal U}$ is isomorphic to ${\mathcal B}$, then $x_i$, is
      essential in $t$ with respect to ${\mathcal U}$ iff $x_i$ is
      essential in $t$ with respect to ${\mathcal B}$. Indeed, if
      $\varphi: {\mathcal U} \rightarrow {\mathcal B}$ is an
      isomorphism and $x_i$ is essential in $t$ with respect to ${\mathcal U}$ then there are elements $a_1, \ldots, a_{i-1}, a_i, a_{i+1},
      \ldots, a_n, ~b_i \not= a_i$ such that $t^{\mathcal U}(a_1,\ldots,a_{i-1},a_i,
      a_{i+1},\ldots,a_n)$ $$\not=t^{\mathcal U}(a_1,\ldots,a_{i-1},b_i,a_{i+1},
      \ldots,a_n).$$ But then\\ \\
       $\varphi(t^{\mathcal U}(a_1,\ldots,
      a_{i-1}, a_i, a_{i+1}, \ldots, a_n))$ \\
      \begin{tabular}{ll}~~~~~~&$=
            t^{\mathcal B}(\varphi(a_1), \varphi(a_2),\ldots, \varphi(a_{i-1}), \varphi(a_i),
      \varphi(a_{i+1}),
      \ldots, \varphi(a_n))$\\ &
     $\not= t^{\mathcal B}(\varphi(a_1), \varphi(a_2), \ldots,
      \varphi(a_{i-1}), \varphi(b_i), \ldots, \varphi(a_n))$
      \end{tabular}\\ \\
       with
      $\varphi(a_i)
      \not= \varphi(b_i)$. Therefore $x_i$ is essential in $t$ with
      respect to the algebra ${\mathcal B}$.
\end{enumerate}
{\it Examples}
\begin{enumerate}
\item The variable $x_i$ is essential in the term $x_i$ with respect to the
      algebra ${\mathcal U}$ of type $\tau$ iff $|A| > 1$.
\item If $t^{\mathcal U}$ is constant then no variable is essential in $t$ with respect
      to $\mathcal U$.
\item If the variable $x_i$ does not occur in the term $t$ then $x_i$ is not
      essential in $t$ with respect to any algebra ${\mathcal U}$ of type
      $\tau$ since $Ess(t, {\mathcal U}) \subseteq var(t)$.
\item Let $\mathcal{SL}$ be a two-element semilattice. Then $x_i$ is
      essential in a term $t$ with respect to $\mathcal{SL}$ iff $x_i \in
      var(t)$ (indeed, terms have the form $t(x_1, \ldots, x_n) = x_{i_1}
      \cdot \ldots \cdot x_{i_n},~ \{i_1, \ldots, i_n \} \subseteq \{1,
      \ldots, n \}$ if the binary operation symbol is written as $\cdot$).
\item Let $\mathcal{BU}$ be a two-element Boolean algebra with conjunction,
      disjunction, and negation ($\land^{\mathcal{BU}}, \lor^{\mathcal{BU}},
      -^{\mathcal{BU}})$ as fundamental operations . Then the
      variable $x_1$ is essential in a term $t$ with respect to
      $\mathcal{BU}$ for the following terms: $x_1, x_1 \lor x_2,~x_2
      \land (x_3 \lor x_1),~(x_1 \land x_2) \lor ({\overline x}_1 \land
      x_3),~(x_1 \land x_2) \lor (x_1 \land {\overline x}_2),~(x_1 \land
      x_2) \lor (x_1 \land x_3) \lor (x_2 \land x_3)$. \\
      But $x_1$ is not essential in $(x_1 \land x_2) \lor ({\overline x}_1
      \land x_2)$ with respect to $\mathcal{BU}$.\\
\end{enumerate}

Variables which are essential in a term $t \in W_\tau(X_n)$ with respect to
an algebra ${\mathcal U} \in Alg(\tau)$ can also be characterized by {\it
evaluation mappings} which we have already introduced in Section \ref{sec2}. We will
use the following more general definition.

\begin{df}\label{3.2}\rm  Let ${\mathcal U} = (A; (f_i^{\mathcal U})_{i \in I})$ be an
algebra of type $\tau$ and let ${\mathcal A}$ be a set with $|{\mathcal A}| = |A|$ and ${\mathcal A} \cap X = \emptyset$. A mapping $h: X_n \cup
{\mathcal A} \rightarrow P_\tau(X_n, {\mathcal A})$ is called {\it
evaluation of the set} $M = \{x_{i_1}, \ldots, x_{i_m} \} \subseteq X_n$ {\it
with the sequence} $C =  ({\overline c}_{i_1}, \ldots, {\overline c}_{i_m} )
\in {\mathcal A}^m$ if
\[ h(x_j) = \left\{ \begin{array}{lll}
   x_j & \mbox{for} & x_j \not\in M \\
\overline{c_j} & \mbox{for} & x_j \in M \end{array}\right.~~\mbox{and}~~h({\overline
   a}) = {\overline a} ~~\mbox{for~every} ~~{\overline a} \in {\mathcal A}\]
\end{df}

\noindent
Clearly, every evaluation of the set $M$ with $C \in {\mathcal A}^m$
is uniquely determined by $M$ and $C$ and can be extended to a uniquely
determined mapping (endomorphism) ${\overline h}: P_\tau(X_n, {\mathcal A}) \rightarrow
P_\tau(X_n, {\mathcal A})$. For the extension there holds ${\overline
h}({\overline h}(p)) = {\overline h}(p)$ for every polynomial $p$ over
${\mathcal A}$, i.e. ${\overline h}$ is idempotent. Moreover we have
$var({\overline h}(p)) = var(p) \setminus M$ ({\it see e.g. \cite{s4,s7,s9}}).

If $M = X_n,~ C \in {\mathcal A}^n$ and if we substitute the
constants ${\overline a} \in {\mathcal A}$ by its corresponding elements
from $A$ we obtain the evaluation mapping of $X$ with $A$ introduced in
Section \ref{sec2}. Note that the result of an evaluation mapping defined by
Definition \ref{3.2} is a polynomial over ${\mathcal A}$.

Now we have

\begin{prop} \label{3.3} A variable $x_i \in X_n$ is essential in the term
$t$ (of type $\tau$) if and only if there exists an evaluation $h$ of the set
$M = X_n \setminus \{ x_i \}$ with some sequence of $C \in {\mathcal A}^{n-1}~(|{\mathcal A}| = A)$ such that the unary polynomial operation
${\overline h}(t)^{\mathcal U}$ takes on at least two values, i.e.
is not constant.   \hfill $\rule{2mm}{2mm}$
\end{prop}

This can also be expressed in the following form: A variable $x_i \in X_n$ is
essential in the term $t$ (of type $\tau$) with respect to the algebra
${\mathcal U}$ (of type $\tau$) iff there exists at least one evaluation
$h$ of $X_n \setminus \{ x_i \}$ with $C \in {\mathcal A}^{n-1}~
(|{\mathcal A}| = |A|)$ such that $x_i \in Ess({\overline h}(t),
{\mathcal U})$.

Another easy consequence of the definition is

\begin{prop}\label{3.4} Let $M \subseteq X_n$ be a nonempty subset of $X_n$
and $x_j \not\in M$. If for every evaluation $h$ of $M$ with $C \in {\mathcal A}^{|M|},~|{\mathcal A}| = |A|$ there holds
$x_j \not\in Ess({\overline h}(t),{\mathcal U})$,
 where ${\mathcal U}$ is an algebra of type
$\tau$,  then $x_j \not\in Ess(t,{\mathcal U})$.
\end{prop}
 ~~~~\hfill $\rule{2mm}{2mm}$

The following Lemma characterizes essential variables with respect to an
algebra ${\mathcal U}$ by non-satisfaction of certain identities in
${\mathcal U}$.

\begin{lemma} \label{3.5} A variable $x_i \in X_n$ is essential in the
$n$-ary term $t$ (of type $\tau$) with respect to an algebra ${\mathcal U}$ of type $\tau$ iff
\[{\mathcal U} \not\models t \approx {\overline h}(t)~,\]
where $h: X_n \rightarrow W_\tau(X_{n+1})$ is a mapping defined by $h(x_i) =
x_{n+1}$ and $h(x_j) = x_j$ for all $j \not= i,~j \in \{ 1, \ldots, n \}$
and where ${\overline h}$ is the extension of $h$, i.e., ${\overline h}: W_\tau(X_n)
\rightarrow W_\tau(X_{n+1})$.
\end{lemma}

\begin{proof} The variable $x_i$ does not occur in the term ${\overline
h}(t)$. Therefore, $x_i$ is not essential for ${\overline h}(t)$ with
respect to ${\mathcal U}$. By Definition \ref{3.1} ${\overline
h}(t)^{\mathcal U}$ does not depend essentially on its $i$-th input $x_i$.

If ${\mathcal U} \models t \approx {\overline h}(t)$ then the induced polynomial
operations are equal: $t^{\mathcal U} = {\overline h}(t)^{\mathcal U}$.
So $t^{\mathcal U}$ does not depend on its $i$-th input $x_i$ and $x_i$ is
not essential for $t$ with respect to ${\mathcal U}$.

If ${\mathcal U} \not \models t \approx {\overline h}(t)$ then there
exists an evaluation $h_1$, of $X_{n+1}$ with $C_1 \in {\mathcal A}^{n+1}, ~
|\mathcal A| = |A|$
such that ${\overline
h}_1(t) \not= \overline{h_1}({\overline h}(t))$. We can assume that $h_1(x_i) \not=
h_1(x_{n+1})$ since $t$ does not depend on $x_{n+1}$ and ${\overline h}(t)$ does not
depend on $x_i$. Now we choose a second evaluation $h_2$ of $X_{n+1}$ with
$C_2 \in {\mathcal A}^{n+1}$ satisfying $h_2(x_i) = h_2(x_{n+1}) =
h_1(x_{n+1}), h_2(x_j) = h_1(x_j)$ if $j \not=
i,n+1$. Clearly, ${\overline h}_2({\overline h}(t)) = {\overline h}_2(t)$ and
${\overline h}_2({\overline h}(t)) = {\overline h}_1({\overline h}(t))$,
therefore ${\overline h}_1(t) \not= {\overline h}_2(t)$. This means, $x_i$
is essential in $t$ with respect  to ${\mathcal U}$.
\end{proof}

Denote by $Id\,{\mathcal U}$ the set of all identities satisfied in
${\mathcal U}$ and by $V({\mathcal U})$ the variety generated by
${\mathcal U}$. It is well-known that $Id\,{\mathcal U}$ is a fully
invariant congruence relation on the absolutely free algebra
$\mathcal{W}_\tau(X_{n+1})$
if the identities in $Id\,{\mathcal U}$ contain at most $n+1$ variables. The
quotient algebra $\mathcal{W}_\tau(X_{n+1})/Id\,{\mathcal U}$
is isomorphic to the
$V({\mathcal U})$-free algebra with $n+1$ generators, i.e., to
$\mathcal{F}_{V({\mathcal U})(X_{n+1})}$.
Then ${\mathcal U} \models t \approx {\overline h}(t)$ is
equivalent to the equality $[t]_{Id\,{\mathcal U}} = [{\overline
h}(t)]_{Id\,{\mathcal U}}$, $i.e.$ the corresponding elements of
$\mathcal{F}_{V({\mathcal U})}(Y)$ where $Y$ is an arbitrary set of free generators
with at least $n+1$ elements agree. But this means:

\begin{cor}\label{3.6} The variable $x_i$ is essential in the $n$-ary term
$t$ with respect to the algebra ${\mathcal U}$ iff $x_i$ is essential in
$t$ with respect to any $V({\mathcal U})$-free algebra with at least $n+1$
free generators.   \hfill $\rule{2mm}{2mm}$
\end{cor}

If the algebra ${\mathcal U}$ is a homomorphic image or a subalgebra of
${\mathcal B}$, then for the sets of identities satisfied in ${\mathcal U}$ and in ${\mathcal B}$, respectively, there holds $Id\,{\mathcal U}
\supseteq Id\,{\mathcal B}$. If ${\mathcal U}$ is a direct power of
${\mathcal B}$ then $Id\,{\mathcal U} = Id\,{\mathcal B}$. Therefore
${\mathcal U} \not \models t \approx {\overline h}(t)$ implies
${\mathcal B} \not \models t \approx {\overline h}(t)$ for the mapping $h$
used in Lemma \ref{3.5}. It follows:

\begin{cor}\label{3.7} If the variable $x_i, 1 \leq i \leq n$, is essential
in the term $t \in W_\tau(X_n)$ with respect to the algebra ${\mathcal U}$
of type $\tau$ then $x_i$ is essential in $t$ with respect to any algebra
${\mathcal B}$ (of type $\tau$) with ${\mathcal U} \in
\Bbb{H}({\mathcal B}),~{\mathcal U} \in \Bbb{S}({\mathcal B})$, or
${\mathcal U} \in \Bbb{P}({\mathcal B})$ where $\Bbb{P}$ is the
operator of forming direct powers. Further $x_i$ is essential for $t$ with
respect to any ${\mathcal B} \in \Bbb{P}({\mathcal U})$.
\hfill $\rule{2mm}{2mm}$
\end{cor}

From Corollary \ref{3.6} we obtain:

\begin{cor}\label{3.8} Let $s,t \in W_\tau(X_n),~n \geq 1$, and assume that
${\mathcal U} \in Alg(\tau)$. If ${\mathcal U} \models s \approx t$ then
$Ess(t, {\mathcal U}) = Ess(s, {\mathcal U})$. \hfill $\rule{2mm}{2mm}$
\end{cor}

Corollary \ref{3.8} suggests the definition of a variable being essential in
a term with respect to a variety.

\begin{df}\label{3.9}\rm  Let $V$ be a variety of type $\tau$ and let $t \in
W_\tau(X_n)$. Then a variable $x_i \in X_n$ is called essential in $t$ with
respect to the variety $V$ if it is essential in $t$ with respect to the
free algebra $\mathcal{F}_V(X)$ with $X = \{ x_1, \ldots, x_n, \ldots \}$
as set of free generators. The set of all variables in $t$ which are
essential with respect to the variety $V$ is denoted by $Ess(t, V)$.
\end{df}

\noindent
{\it Example } The class $Alg(\tau)$ is the biggest variety of algebras of
type $\tau$. The free algebra $\mathcal{W}_\tau(X)$ with respect to
$Alg(\tau)$ is the algebra of all terms of type $\tau$ where the operations
$f_i^{\mathcal{W}_\tau(X)}$ are defined by $(t_1, \ldots, t_{n_i})
\longmapsto f_i^{\mathcal{W}_\tau(X)} (t_1, \ldots, t_{n_i}) := f_i(t_1,
\ldots, t_{n_i})$ ({\it see e.g.\cite{s8}}).\\

A variable $x_i \in X_n$ is essential in the term $t$ with
respect to the variety $Alg(\tau)$ iff $x_i \in var(t)$. Indeed,
if $x_i \in Ess(t, Alg(\tau))$ then $x_i \in var(t)$. Conversely, if $x_i \in
var(t)$ then $t \not= {\overline h}(t)$ for the mapping $h: X_n \rightarrow
W_\tau(X_{n+1})$ defined by $h(x_i) = x_{n+1}$ and $h(x_j) = x_j$ for all $j
\not= i$. Therefore, $Alg(\tau) \not \models t \approx {\overline h}(t)$ and
$x_i \in Ess(t, Alg(\tau))$.

\begin{prop}\label{3.10} If $x_i \in X_n$ is essential in the $n$-ary term
$t$ of type $\tau$ with respect to the variety $V$ of type $\tau$ and if $W
\supseteq V$, i.e. if $V$ is a subvariety of $W$ then $x_i$ is essential in
$t$ with respect to $W$.
\end{prop}

\begin{proof} $W \supseteq V$ implies $Id\,W \subseteq Id\,V$. If $x_i \in
Ess(t,V)$ then $V \not \models t \approx {\overline h}(t)$ for the mapping
defined in Lemma \ref{3.5}. But then $W \not \models t \approx {\overline h}(t)$
and $x_i \in Ess(t,W)$.  \end{proof}

\noindent
{\it Example} We consider the lattice of all varieties of semigroups. An
equation $s \approx t$ is called regular if $var(s) = var(t)$. A variety is
regular if it is the model class of a set of regular equations. It is
well-known that a variety $V$ of semigroups is regular iff it contains the
variety $SL$ of semilattices and that the variety of semilattices is an
atom in the lattice of all varieties of semigroups (\cite{s6}). By
 Example 4 after Definition \ref{3.1} we have $x_i \in Ess(t, SL)$ iff $x_i \in var(t)$.
Now Proposition \ref{3.10} shows that for an arbitrary regular variety
$V$ of semigroups we have $x_i \in Ess(t, V)$ iff $x_i \in var(t)$.\\

A term $t \in W_\tau(X)$ is called a subterm of a term $s \in W_\tau(X)$
with respect to the algebra ${\mathcal U}$ of type $\tau$ and we write $t
\prec s$ if there exists an evaluation $h$ of some subset $M \subset var(s)$
with a sequence $C \in {\mathcal A}^{|M|},~(|{\mathcal A}| = |A|)$ such
that ${\mathcal U} \models t \approx {\overline h}(s)$. By $Sub(t,
{\mathcal U})$ we denote the set of all subterms of a term $t$ of type
$\tau$ with respect to the algebra ${\mathcal U}$.

Clearly, if $t \prec s$ then $Ess(t, {\mathcal U}) \subseteq Ess(s,
{\mathcal U})$.

The concept of a {\it separable set} of variables in a function defined e.g.
in \cite{s2} can also be extended to terms.

\begin{df}\label{3.11}\rm Let ${\mathcal U}$ be an algebra of type $\tau$ and
let ${\mathcal A}$ be a set of constant symbols (for elements of $A$) with
$|A| = |{\mathcal A}|$. A set $M$ of essential variables in the term $t \in
W_\tau(X)$ with respect to the algebra $\mathcal U$
is called {\it separable in $t$
with respect to $\mathcal U$} if there exists at least one evaluation ${\overline h}$ of $X_n
\setminus M$ with $C \in {\mathcal A}^{n-|M|}$ such that $M = Ess({\overline
h}(t),{\mathcal U})$.\\

By $Sep(t,{\mathcal U})$ we denote the set of all separable sets in $t$
with respect to ${\mathcal U}$.
\end{df}

\noindent
{\it Example}\   Consider the term $t_1= x_1 x_2 + x_3$ with respect to the
two-element Boolean ring $\mathcal{BR} = (\{0,1\}; + , \cdot)$. Then $M =
\{ x_1, x_2 \}$ is separable in $t$ with respect to $\mathcal{BR}$ since
$h: \{ x_3 \} \rightarrow {\overline 0}$ gives ${\overline h}(t_1) = x_1
x_2$ and $Ess({\overline h}(t_1), \mathcal{BR}) = \{ x_1, x_2 \}$.\\

For separable sets of variables with respect to an algebra ${\mathcal U}$
we obtain results which are similar to some results on essential variables.
We add some facts about separable sets. The proofs are straightforward and left
to the reader.

\begin{theo}\label{3.12} Let $s, t$ be terms of type $\tau$ and let
${\mathcal U}$ be an algebra of type $\tau$. \\[-7mm]
\begin{enumerate}
\item[(i)] If $s$ and $t$ are two terms with ${\mathcal U} \models s \approx t$
      then $Sep(s, {\mathcal U}) = Sep(t, {\mathcal U})$.
\item[(ii)] If $s \prec t$ then $Sep(s, {\mathcal U}) \subseteq Sep(t,
      {\mathcal U})$
\item[(iii)] A set $M = \{x_{i_1}, \ldots, x_{i_m} \} \subseteq X_n$ of essential
      variables in the term $t$ with respect to the algebra ${\mathcal U}$
      is separable in $t$ with respect to ${\mathcal U}$ iff $M$ is
      separable in $t$ with respect to any free algebra in $V({\mathcal U})$ with at least $n+1$ free generators. \hfill
                      $\Box$
\end{enumerate}
\end{theo}

The last proposition motivates the following definition:

\begin{df}\label{3.13}\rm Let $t$ be an $n$-ary term of type $\tau$ and let
$V$ be a variety of type $\tau$. A set $M \subseteq X_n$ is called separable
in $t$ with respect to $V$ if $M$ is separable in $t$ with respect to the
free algebra $\mathcal{F}_V(X_n)$. Let $Sep(t,V)$ be the set of all
separable sets with respect to $V$.
\end{df}

\section{Complexity of terms, polynomials, and algebras}

\noindent
Terms are useful tools in Theoretical Computer Science. They can be used as
models for different structures in logic programming. Every term can be
described by a graph with operation symbols or variables as vertices and with
subterms as edges. Such graphs are also called {\it semantic trees} of the
corresponding terms.\\

\noindent
{\it Example:} Figure \ref{f1} shows the semantic trees of the terms $t_1 = x_1x_2 + x_3$
and $t_2 = x_1x_3 + x_2\bar{x_3}$ of a language containing two binary
operation symbols $\cdot$ and $+$ and a unary symbol $-$.\\
\begin{figure}
\begin{center}
\unitlength=0.50mm
\special{em:linewidth 0.4pt}
\linethickness{0.4pt}
\begin{picture}(115.00,71.00)
\put(22.00,68.00){\vector(-1,-2){9.83}}
\put(11.75,48.00){\vector(-1,-2){9.67}}
\put(22.00,68.00){\vector(1,-2){9.67}}
\put(11.75,48.00){\vector(1,-2){9.75}}
\put(22.00,68.00){\circle*{1.00}}
\put(11.75,48.00){\circle*{1.00}}
\put(23.00,74.00){\makebox(0,0)[cc]{$t_1$}}
\put(17.00,67.00){\makebox(0,0)[cc]{$+$}}
\put(9.00,48.00){\makebox(0,0)[cc]{$.$}}
\put(32.50,45.50){\makebox(0,0)[cc]{$x_3$}}
\put(23.00,24.00){\makebox(0,0)[cc]{$x_2$}}
\put(3.00,24.00){\makebox(0,0)[cc]{$x_1$}}
\put(32.00,48.00){\circle*{1.00}}
\put(21.83,28.50){\circle*{1.00}}
\put(2.00,28.00){\circle*{1.00}}
\put(61.75,48.00){\vector(-1,-2){9.67}}
\put(61.75,48.00){\vector(1,-2){9.75}}
\put(61.75,48.00){\circle*{1.00}}
\put(59.00,49.00){\makebox(0,0)[cc]{$.$}}
\put(72.55,24.00){\makebox(0,0)[cc]{$x_3$}}
\put(53.00,24.00){\makebox(0,0)[cc]{$x_1$}}
\put(71.55,28.00){\circle*{1.00}}
\put(52.00,28.00){\circle*{1.00}}
\put(102.33,48.00){\vector(-1,-2){9.75}}
\put(102.00,48.00){\vector(1,-2){9.75}}
\put(102.00,48.00){\circle*{1.00}}
\put(105.00,49.00){\makebox(0,0)[cc]{$.$}}
\put(115.00,26.00){\makebox(0,0)[cc]{$^-$}}
\put(93.25,24.00){\makebox(0,0)[cc]{$x_2$}}
\put(92.25,28.25){\circle*{1.00}}
\put(112.00,28.00){\circle*{1.00}}
\put(113.00,11.00){\makebox(0,0)[cc]{$x_3$}}
\put(112.00,14.50){\circle*{1.00}}
\put(82.00,68.00){\vector(-1,-1){20.00}}
\put(82.00,68.00){\vector(1,-1){20.00}}
\put(82.00,68.00){\circle*{1.00}}
\put(87.00,68.00){\makebox(0,0)[cc]{$+$}}
\put(83.00,74.00){\makebox(0,0)[cc]{$t_2$}}
\put(112.00,28.00){\vector(0,-1){13.00}}
\end{picture}
\end{center}
\caption{}\label{f1}
\end{figure}

This example shows that the graph of $t_2$ is more complex than
the graph of $t_1$ since the number of edges and vertices in the
graph of $t_2$  is greater than in $t_1$. There are different kinds of
concepts of {\it complexity} of terms. Roughly spoken, there are two
classes of term complexity, i.e. of functions from $W_{\tau}(X) \to \Bbb N$.
The first one is based on a linguistic point of view and counts the number
of variables or of operation symbols occurring  in the term.
Let us give two examples. If $\ell_i$ denotes the number of occurrences of
the variable $x_i$ in the $n$-ary term $t$ then the first complexity measure
is defined as
\[Cp^{(1)}(t) = \sum \limits_{i=1}^{n} \ell_i .\]
The second one counts the operation symbols occurring in the term and is
defined inductively by
\[Cp^{(2)}(x_i) := 0,~ Cp^{(2)}(f_i(t_1, \ldots, t_{n_i})):= \sum \limits_{j=1}^{n_i}
Cp^{(2)}(t_j) + 1.\]
For our terms $t_1,~ t_2$ we obtain $Cp^{(1)}(t_1) = 3,~ Cp^{(1)}(t_2) = 4$~
and~
\[Cp^{(2)}(t_1) = 2,~  Cp^{(2)}(t_2) = 4.\]
The second class is based on the consideration of a term as a word in a given
language together with its ``meaning '' in an algebra.
In the theory of {\it switching circuits} one uses so-called
{\it binary decision diagrams} (bdd) to measure the complexity of terms.
In this case the underlying algebra is the two-element Boolean algebra
$\mathcal{BU} =(\{0,1\}; +, \cdot,-)$ of type (2,2,1) (or another two-element
algebra). One starts to give the variable with the greatest index the value $0$ or $1$,
respectively, and continues with the arising subterms in the same way.
In our examples we obtain the  diagrams given on the Figure \ref{f2}.\\
\begin{figure}
\begin{center}
\unitlength=1.00mm
\linethickness{0.4pt}
\begin{picture}(93.50,55.00)
\put(36.00,20.00){\vector(1,-2){5.00}}
\put(36.00,20.00){\vector(-1,-2){5.00}}
\put(31.00,30.00){\vector(1,-2){5.00}}
\put(31.00,30.00){\vector(-1,-2){5.00}}
\put(9.00,30.00){\vector(1,-2){5.00}}
\put(9.00,30.00){\vector(-1,-2){5.00}}
\put(14.00,20.00){\vector(1,-2){5.00}}
\put(14.00,20.00){\vector(-1,-2){5.00}}
\put(20.00,52.00){\vector(1,-2){11.00}}
\put(20.00,52.00){\vector(-1,-2){11.00}}
\put(20.00,56.00){\makebox(0,0)[cc]{$t_1=x_1x_2+x_3$}}
\put(37.00,30.00){\makebox(0,0)[cc]{$\overline{x_1x_2}$}}
\put(27.50,42.00){\makebox(0,0)[cc]{1}}
\put(12.50,42.00){\makebox(0,0)[cc]{0}}
\put(15.00,30.00){\makebox(0,0)[cc]{$x_1x_2$}}
\put(4.50,25.00){\makebox(0,0)[cc]{0}}
\put(13.50,25.00){\makebox(0,0)[cc]{1}}
\put(26.50,25.00){\makebox(0,0)[cc]{0}}
\put(35.50,25.00){\makebox(0,0)[cc]{1}}
\put(18.00,20.00){\makebox(0,0)[cc]{$x_1$}}
\put(4.00,17.00){\makebox(0,0)[cc]{0}}
\put(26.00,17.00){\makebox(0,0)[cc]{1}}
\put(40.00,20.00){\makebox(0,0)[cc]{$\overline{x_1}$}}
\put(40.50,15.00){\makebox(0,0)[cc]{1}}
\put(31.50,15.00){\makebox(0,0)[cc]{0}}
\put(18.50,15.00){\makebox(0,0)[cc]{1}}
\put(9.50,15.00){\makebox(0,0)[cc]{0}}
\put(9.00,7.00){\makebox(0,0)[cc]{0}}
\put(19.00,7.00){\makebox(0,0)[cc]{1}}
\put(31.00,7.00){\makebox(0,0)[cc]{1}}
\put(41.00,7.00){\makebox(0,0)[cc]{0}}
\put(36.00,20.00){\circle*{1.00}}
\put(41.00,10.00){\circle*{1.00}}
\put(31.00,10.00){\circle*{1.00}}
\put(19.00,10.00){\circle*{1.00}}
\put(9.00,10.00){\circle*{1.00}}
\put(4.00,20.00){\circle*{1.00}}
\put(14.00,20.00){\circle*{1.00}}
\put(26.00,20.00){\circle*{1.00}}
\put(31.00,30.00){\circle*{1.00}}
\put(9.00,30.00){\circle*{1.00}}
\put(20.00,52.00){\circle*{1.00}}
\put(88.00,21.00){\vector(1,-2){5.00}}
\put(88.00,21.00){\vector(-1,-2){5.00}}
\put(66.00,21.00){\vector(-1,-2){5.00}}
\put(77.00,43.00){\vector(1,-2){11.00}}
\put(77.00,43.00){\vector(-1,-2){11.00}}
\put(77.00,47.00){\makebox(0,0)[cc]{$t_2=x_1x_3+x_2\overline{x_3}$}}
\put(93.00,21.00){\makebox(0,0)[cc]{$x_1$}}
\put(84.50,33.00){\makebox(0,0)[cc]{1}}
\put(69.50,33.00){\makebox(0,0)[cc]{0}}
\put(61.50,16.00){\makebox(0,0)[cc]{0}}
\put(81.50,16.00){\makebox(0,0)[cc]{0}}
\put(92.50,16.00){\makebox(0,0)[cc]{1}}
\put(93.00,11.00){\circle*{1.00}}
\put(61.00,11.00){\circle*{1.00}}
\put(71.00,11.00){\circle*{1.00}}
\put(83.00,11.00){\circle*{1.00}}
\put(88.00,21.00){\circle*{1.00}}
\put(66.00,21.00){\circle*{1.00}}
\put(77.00,43.00){\circle*{1.00}}
\put(61.00,8.00){\makebox(0,0)[cc]{0}}
\put(83.00,8.00){\makebox(0,0)[cc]{0}}
\put(93.00,8.00){\makebox(0,0)[cc]{1}}
\put(71.00,21.00){\makebox(0,0)[cc]{$x_2$}}
\put(71.00,16.00){\makebox(0,0)[rc]{1}}
\put(71.00,8.00){\makebox(0,0)[cc]{1}}
\put(66.00,21.00){\vector(1,-2){5.00}}
\end{picture}
\end{center}
 \caption{}\label{f2}
\end{figure}
The bdd of $t_1$ looks more complex. This reflects the fact that in the first case
there are more evaluation mappings.
 If one starts with another variable instead of $x_3$ one gets different bdd's.
 Figure \ref{f3} gives the bdd's of $t_1$ and $t_2$ if we replace the variables following the
 order $x_2, x_3,x_1$ by $0$ and $1$. The reason is that $\{x_1,x_2\}$ is not eparable in $t_2$.\\
\begin{figure}
\begin{center}
\unitlength=1.00mm
\linethickness{0.4pt}
\begin{picture}(97.00,55.00)
\put(36.00,20.00){\vector(1,-2){5.00}}
\put(36.00,20.00){\vector(-1,-2){5.00}}
\put(31.00,30.00){\vector(1,-2){5.00}}
\put(9.00,30.00){\vector(1,-2){5.00}}
\put(9.00,30.00){\vector(-1,-2){5.00}}
\put(20.00,52.00){\vector(1,-2){11.00}}
\put(20.00,52.00){\vector(-1,-2){11.00}}
\put(20.00,55.00){\makebox(0,0)[cc]{$t_1=x_1x_2+x_3$}}
\put(40.00,30.00){\makebox(0,0)[cc]{$x_1+x_3$}}
\put(27.50,42.00){\makebox(0,0)[cc]{1}}
\put(12.50,42.00){\makebox(0,0)[cc]{0}}
\put(14.00,30.00){\makebox(0,0)[cc]{$x_3$}}
\put(4.50,25.00){\makebox(0,0)[cc]{0}}
\put(13.50,25.00){\makebox(0,0)[cc]{1}}
\put(35.50,25.00){\makebox(0,0)[cc]{1}}
\put(4.00,17.00){\makebox(0,0)[cc]{0}}
\put(43.00,20.00){\makebox(0,0)[cc]{$x_1+1$}}
\put(40.50,15.00){\makebox(0,0)[cc]{1}}
\put(31.50,15.00){\makebox(0,0)[cc]{0}}
\put(31.00,7.00){\makebox(0,0)[cc]{1}}
\put(41.00,7.00){\makebox(0,0)[cc]{0}}
\put(36.00,20.00){\circle*{1.00}}
\put(41.00,10.00){\circle*{1.00}}
\put(31.00,10.00){\circle*{1.00}}
\put(4.00,20.00){\circle*{1.00}}
\put(14.00,20.00){\circle*{1.00}}
\put(31.00,30.00){\circle*{1.00}}
\put(9.00,30.00){\circle*{1.00}}
\put(20.00,52.00){\circle*{1.00}}
\put(88.00,21.00){\vector(1,-2){5.00}}
\put(88.00,21.00){\vector(-1,-2){5.00}}
\put(66.00,21.00){\vector(1,-2){5.00}}
\put(66.00,21.00){\vector(-1,-2){5.00}}
\put(77.00,43.00){\vector(1,-2){11.00}}
\put(77.00,43.00){\vector(-1,-2){11.00}}
\put(77.00,46.00){\makebox(0,0)[cc]{$t_2=x_1x_3+x_2\overline{x_3}$}}
\put(99.00,21.00){\makebox(0,0)[cc]{$x_1x_3+\overline{x_3}$}}
\put(75.00,11.00){\makebox(0,0)[cc]{$x_1$}}
\put(84.50,33.00){\makebox(0,0)[cc]{1}}
\put(69.50,33.00){\makebox(0,0)[cc]{0}}
\put(73.00,21.00){\makebox(0,0)[cc]{$x_1x_3$}}
\put(61.50,16.00){\makebox(0,0)[cc]{0}}
\put(70.50,16.00){\makebox(0,0)[cc]{1}}
\put(83.50,16.00){\makebox(0,0)[cc]{0}}
\put(92.50,16.00){\makebox(0,0)[cc]{1}}
\put(93.00,11.00){\circle*{1.00}}
\put(61.00,11.00){\circle*{1.00}}
\put(71.00,11.00){\circle*{1.00}}
\put(83.00,11.00){\circle*{1.00}}
\put(88.00,21.00){\circle*{1.00}}
\put(66.00,21.00){\circle*{1.00}}
\put(77.00,43.00){\circle*{1.00}}
\put(61.00,7.00){\makebox(0,0)[cc]{0}}
\put(83.00,8.00){\makebox(0,0)[cc]{0}}
\put(93.00,8.00){\makebox(0,0)[cc]{1}}
\put(31.00,30.00){\vector(-1,-1){10.00}}
\put(21.00,20.00){\vector(1,-2){5.00}}
\put(21.00,20.00){\vector(-1,-2){5.00}}
\put(25.50,15.00){\makebox(0,0)[cc]{1}}
\put(16.50,15.00){\makebox(0,0)[cc]{0}}
\put(16.00,7.00){\makebox(0,0)[cc]{0}}
\put(26.00,7.00){\makebox(0,0)[cc]{1}}
\put(21.00,20.00){\circle*{1.00}}
\put(26.00,10.00){\circle*{1.00}}
\put(16.00,10.00){\circle*{1.00}}
\put(14.00,17.00){\makebox(0,0)[cc]{1}}
\put(24.00,25.00){\makebox(0,0)[cc]{0}}
\put(26.00,20.00){\makebox(0,0)[cc]{$x_1$}}
\put(71.00,11.00){\vector(-1,-2){5.00}}
\put(71.00,11.00){\vector(1,-2){5.00}}
\put(66.00,1.00){\circle*{1.00}}
\put(76.00,1.00){\circle*{1.00}}
\put(66.00,8.00){\makebox(0,0)[cc]{0}}
\put(76.00,6.00){\makebox(0,0)[cc]{1}}
\put(66.00,-1.50){\makebox(0,0)[cc]{0}}
\put(76.00,-1.50){\makebox(0,0)[cc]{1}}
\end{picture}
\end{center}
 \caption{}\label{f3}
\end{figure}

Now we want to count all evaluations of
$X_n \setminus M$
for arbitrary (not only one-element) subsets $M$ of $X_n$ with elements
of an arbitrary algebra $\mathcal U$ of type $\tau$. While $Cp^{(1)}$ and $Cp^{(2)}$ are complexities of terms,
$Cp^{(3)}$ as the complexity of a function, is connected with the computational complexity.

\begin{df} \label{4.1}\rm
Let $t \in W_{\tau}(X_n)$ be a term and let $\emptyset \not= M \subseteq X_n$
be a subset. The number $Cp^{(3)}(t,M,\mathcal U)$ of all those evaluations
of $X_n \setminus M$ with $C \in {\mathcal A}^{n - |M|}, |A| = |\mathcal A|$, for which $M =
Ess(\overline{h}(t), \mathcal U)$ is called {\it complexity of}\   $t$ {\it with
respect to $M$ in $\mathcal U$.}
Further we define the complexity of\  $t$ in $\mathcal U$ by
\[Cp^{(3)}(t,\mathcal U) := \sum \limits_{\emptyset \subset M \subseteq X_n}
Cp^{(3)}(t,M, \mathcal U).\]
\end{df}

\noindent
{\it Example:} We consider again the terms $t_1, t_2$ in the two-element Boolean algebra
$\mathcal{BU} =(\{0,1\}; +,\cdot, -)$ and obtain\\
$Cp^{(3)}(t_2, \{x_3\}, \mathcal{BU}) = 2,
 Cp^{(3)}(t_2, \{x_2\}, \mathcal{BU}) = 2,
 Cp^{(3)}(t_2, \{x_1\}, \mathcal{BU}) = 2,$
$Cp^{(3)}(t_2, \{x_1,x_2\}, \mathcal{BU}) = 2,
 Cp^{(3)}(t_2, \{x_1,x_3\}, \mathcal{BU}) = 2,\\
 Cp^{(3)}(t_2, \{x_2,x_3\}, \mathcal{BU}) = 0,$
 $Cp^{(3)}(t_2, \{x_1,x_2,x_3\}, \mathcal{BU}) = 1$\\
  and therefore
$ Cp^{(3)}(t_2, \mathcal{BU}) = 11$.\\

In a similar way we calculate $Cp^{(3)}(t_1, \mathcal{BU}) = 13$.

Since we defined $Cp^{(3)}(t, \mathcal U)$ using evaluation mappings we can
expect that for terms $t,t'$ which form an identity in $\mathcal U$ we
get the same complexity.

\begin{lemma} \label{4.2}
If $\mathcal U \models t \approx t'$ then for each
$\emptyset \subset M \subseteq X_n$ we have
\[Cp^{(3)}(t,M,\mathcal U) =  Cp^{(3)}(t',M,\mathcal U)\]
and thus also \[Cp^{(3)}(t,\mathcal U) =  Cp^{(3)}(t',\mathcal U).\]
\end{lemma}
\begin{proof} If $\mathcal U \models t \approx t'$ then
$Ess(t,\mathcal U) = Ess(t', \mathcal U)$ by Corollary \ref{3.8} of Lemma \ref{3.5}. Let $h$ be an evaluation
of $X_n\setminus M$ with $C \in {\mathcal A}^{n-|M|}, (|\mathcal A| = |A|)$ satisfying
$M = Ess(\overline{h}(t), \mathcal U)$. From $\mathcal U \models t \approx t'$
for every evaluation of $X_n\setminus M$ with elements from $\mathcal A$
we obtain  polynomials $\overline{h}(t),\overline{h}(t')$ which give equal
polynomial operations $\overline{h}(t)^{\mathcal U} =
\overline{h}(t')^{\mathcal U}$, $i.e.$ which form a polynomial identity. But this means
$Ess(\overline{h}(t),\mathcal U) = Ess(\overline{h}(t'), \mathcal U)$.
Therefore, $h$ is an evaluation of $X_n\setminus M$ with
$C \in {\mathcal A}^{n -|M|}$ satisfying
$M = Ess(\overline{h}(t),\mathcal U)$ iff $h$ is an evaluation of $X_n\setminus M$
with $\mathcal A$ satisfying
$M = Ess(\overline{h}(t'),\mathcal U)$. \end{proof}

This motivates to consider only the elements of the free
algebra $F_{V(\mathcal U)}(X_n)$, freely generated by $X_n$.
Let $t^*$ be the equivalence class with respect to the equivalence
relation $Id \mathcal U$ containing the term $t$. Then we define

\begin{df} \label{4.3}\rm The sum
\[\sum\limits_{t^* \in F_{V(\mathcal U)}(X_n)} Cp^{(3)}(t^*, \mathcal U)\]
is called {\it $n$-complexity of $\mathcal U$}.
\end{df}

\noindent
{\it Example} A finite algebra $\mathcal U = (A; (f_i^{\mathcal U})_{i \in I})$
is called {\it primal} if $T(\mathcal U) = O_A$, $i.e.$ if every operation
defined on $A$ is a term operation of $\mathcal U$. We consider the variety
$V(\mathcal U)$ generated by a primal algebra with $|A| = k \geq 2$.
Let $\mathcal{F}_{V(\mathcal U) (X_n)}$ be the free algebra in
$V(\mathcal U)$ generated by an $n$-element alphabet $X_n$.\\

Then there are exactly $k^{k^n}$ different equivalence
classes with respect to the equivalence relation
$Id{\mathcal U}$, $i.e.$ $|\mathcal F_{V(\mathcal U)}(X_n)| = k^{k^n}$.
The algorithm for calculating  the complexity of
$\mathcal U$ can be described in the following way:
Let us set $p=k^{k^n}-1$ and $s=2^n.$

\vspace{1cm}

\noindent
begin \\
1\hspace{2cm} Coding all
term operations from $\mathcal U$
by the \\
2\hspace{3cm}  integers from the set $\{0,1,\ldots, p\}$;       \\
3\hspace{2cm}  CPA:=0;                                                  \\
4\hspace{2cm}  for $i:=0$ to $p$ do                                           \\
5\hspace{3cm}  begin \\
6\hspace{4cm}  $CP:=0$; \\
7\hspace{4cm}  Coding all subsets of $X_n$ by the integers \\
8\hspace{5cm}  from the set $\{1,2,\ldots,s\}$;\\
9\hspace{4cm}  for $j:=1$ to $s$ do                    \\
10\hspace{5cm}  begin                             \\
11\hspace{5.5cm}  $r:=k^{|Decode(j)|-1}$;\\
12\hspace{5.5cm}  Coding all evaluations of the \\
13\hspace{6cm}  set $X_n\setminus Decode(j)$\\
14\hspace{6cm}  by the integers from $\{0,\ldots,r\}$;\\
15\hspace{5.5cm}  for $m:=0$ to $r$ do \\
16\hspace{5.5cm}  if $Decode(j)=Ess(Decode(i,m),\mathcal U)$ \\
17\hspace{5.5cm}  then $CP:=CP+1$; \\
18\hspace{5.5cm}  Print $CP$; \\
19\hspace{5cm}  end;     \\
20\hspace{5cm}  $CPA:=CPA+CP$;  \\
21\hspace{3cm}  end; \\
22\hspace{3cm}  Print $CPA$ as complexity of $\mathcal U$; \\
end.\\

Here $Decode(j)$ is the subset of $X_n$ which is decoded by $j$, and
$Decode(i,m)$ is the image of the polynomial coded by $i$ under the
evaluation coded by $m$.

This algorithm is constructed by full exhaustion of all cases and therefore
it works not so quickly, but for small values of $k$ and $n$ it can be used.
In the case
$k=2$ and $n=3$ we obtained the following results.
All $n$-ary functions (unequivalent classes of terms over $V(\mathcal U))$
are classified with respect to the complexity.
There exist 2 terms with
complexity 19, 16 with complexity 16, 40  with complexity 13,
72  with complexity 12, 24  with complexity 11,
6  with complexity 10, 48  with complexity 9,
24  with complexity 8, 16  with complexity 7,
6  with complexity 4 and
2 termss with complexity 0.

Consequently, the $3$-complexity of the two-element Boolean algebra
$\mathcal{BU}$ is:\\

\begin{tabular}{ll}
$Cp^{(3)}(\mathcal U)=$&$2\cdot19+16\cdot16+40\cdot13+72\cdot12+24\cdot11$\\ & $+
6\cdot10+48\cdot9+24\cdot8+16\cdot7+6\cdot4+2\cdot0=$\\ &
$=38+256+520+864+264+60+432$\\ &$+192+112+24+0=2762.$\\& \\
\end{tabular}

We can also define the $n$-complexity of a variety.

\begin{df}\label{4.4}\rm
Let $V$ be a variety of type $\tau$ and let $\mathcal{F}_V(X)$ be the $V$-free algebra
generated by the $n$-element alphabet $X_n$. Then
\[\sum \limits_{t^* \in F_V(X_n)} Cp^{(3)}(t^*,\mathcal{F}_V(X))\] is
called {\it $n$-complexity of $V$}.
\end{df}

Lemma \ref{4.2} gave an example which showed that different terms can have the
same complexity. The following results give examples for equal complexities
if we have different $M's$ or different algebras.

Let $\sigma$ be a permutation of $X_n$ and let $f_\sigma: W_{\tau}(X_n) \to
W_{\tau}(X_n)$ be the extension of $\sigma$ to terms. Then $f_{\sigma}$
has the following interesting property.
For each subset $\emptyset \subset M \subseteq X_n$, for each term
$t \in W_{\tau}(X_n)$ and for each algebra $\mathcal U$ of type $\tau$
we have \[Cp^{(3)}(t,M,\mathcal U) =
Cp^{(3)}(f_{\sigma}(t), \sigma(M), \mathcal U).\]
For the proof we have to show that for every evaluation $h: X_n \setminus M
\to \mathcal A, |A| =|\mathcal A|$, satisfying  $M= Ess(\overline{h}(t), \mathcal U)$
we get an evaluation $h': X_n\setminus\sigma(M) \to \mathcal A$ satisfying
$\sigma(M) = Ess(\overline{h'}(f_{\sigma}(t)), \sigma(M), \mathcal U)$ and
conversely. Clearly, $X_n\setminus\sigma(M) =\{\sigma(x_i)\ |\ x_i \not \in M\}$.
We define $h'$ by $h'(\sigma(x_i)) := h(x_i)$. It is easy to see that $h'$
has the desired  property.

Let $p$ be a polynomial of type $\tau$ over $\mathcal A$ and let $\mathcal U$
be an algebra of type $\tau$ with $|\overline{A}| = |A|$.
By $V_p^{\mathcal U}$ we denote the set of all values of the induced
polynomial operation $p^{\mathcal U}$, $i.e.$  $$V_p^{\mathcal U}:=\{a \in a\ |\
\exists h: X \cup \mathcal A \to A\  ~(\overline h(p) = a)\}.$$
Let $g: \mathcal A \to \mathcal A$ be a mapping on the set of constant symbols.
This mapping can be uniquely extended to a mapping
$f_g: P_{\tau}(X_n, \mathcal A) \to P_{\tau}(X_n, \mathcal A)$
with $f_g(x_i) = x_i$ for all $i \in \{x_1, \ldots, x_n\}$. Then we have

\begin{prop} \label{4.5}
Let $g: \mathcal A \to \mathcal A$ be a mapping and let
$p \in P_{\tau}(X_n, {\mathcal A})$.
If $|V_p^{\mathcal U}| = |g(V_p^{\mathcal U})|$
then for every
$\emptyset \subset M \subseteq X_n$ we have
\[Cp^{(3)}(p, M, \mathcal U) = Cp^{(3)}(f_g(p), M, \mathcal U)\]
and therefore also
\[Cp^{(3)}(p, \mathcal U) = Cp^{(3)}(f_g(p), \mathcal U).\]
\end{prop}
\begin{proof} For the proof we have to show that for every evaluation
$h:X_n \setminus M \to \mathcal A$, satisfying
$M = Ess(\overline h(p), M, \mathcal U)$ we have an evaluation
$h':X_n \setminus M \to \mathcal A$, satisfying
$M = Ess(\overline{h'}(f_g(p)), M, \mathcal U)$.
Because of $|V_p^{\mathcal U}| = |g(V_p^{\mathcal U})|$
and the finiteness, the restriction of $g$ to $V_p^{\mathcal U}$
is a one-to-one mapping.
Therefore, $f_{g^{-1}}(p)$ exists and satisfies $f_{g^{-1}}(f_g(p)) = p$.
If we define $\overline{h'}:= \overline h \circ f_{g^{-1}}$ then we have
\[\overline{h'}(f_g(p)) = (\overline h \circ f_{g^{-1}})(f_g(p)) =
\overline h(p).\]
Therefore, \[Ess(\overline{h'}(f_g(p)), M, \mathcal U) =
Ess(\overline{h}(p), M, \mathcal U) = M.\]
Clearly, there is a one-to-one correspondence between evaluations
$h$ and $h'$. This proves Proposition \ref{4.5}. \end{proof}

\begin{theo} \label{4.6}
Let $\mathcal U$ and $\mathcal B$ be isomorphic algebras of type $\tau$.
If $t \in W_{\tau}(X)$ and $\emptyset \subset M \subseteq X_n$ then
\[Cp^{(3)}(t,M,{\mathcal U}) = Cp^{(3)}(t,M,\mathcal B).\]
\end{theo}
\begin{proof} Let $\varphi:\mathcal U \to \mathcal B$ be an isomorphism and
let $h: X_n \setminus M  \to \mathcal A$ be an evaluation satisfying
$M = Ess(\overline{h}(t), \mathcal U)$. Consider a set $\overline B$ of constant symbols
for the elements of $B$. Then $\overline{\varphi}: \mathcal A \to \overline B$
is a one-to-one mapping.
Consider the evaluation $h'$ of $X_n \setminus M$ with $\overline B$
defined by $h' = \overline{\varphi} \circ h$.
Since $\overline{\varphi}$ is one-to-one we have
$$M = Ess(\overline{h}(t), \mathcal U) = Ess(\overline{\overline{\varphi}
\circ h}(t), \mathcal B)$$ where $\overline{\overline{\varphi} \circ h}$
is the extension of $\overline{\varphi} \circ h$.
To every $h$ we obtain exactly one $\overline{\varphi} \circ h$.
If conversely $h^*: X_n \setminus M \to \overline B$ is an evaluation then by
$(h^{*})':= {\overline \varphi}^{-1} \circ h$ we obtain an evaluation
$(h^{*})': X_n \setminus M \to \mathcal A$ which is uniquely determined
and satisfies
$$M = Ess(\overline{h^*}(t), \mathcal B) =   Ess(\overline{(h^{*})'}(t),
\mathcal U).$$ Therefore, the cardinalities of the sets of all evaluations of
$X_n \setminus M$ with $\mathcal A$ and $\overline B$ with the property that the
set of essential variables is M, are equal and this means,
\[Cp^{(3)}(t,M,\mathcal U) = Cp^{(3)}(t,M,\mathcal B).\] \end{proof}

\vspace*{3mm}

\end{document}